# A generalization of a result of Hille on analytic implicit functions

A GENERALIZATION OF A RESULT OF HILLE

ON ANALYTIC IMPLICIT FUNCTIONS

by

J.F. Barrett

Contents        Page

Abstract

Hille in 1959 gave a stronger form for the classical implicit function theorem of Cauchy for equations of the type

$$y = \psi(x,y) \qquad x,y \in \mathbb{R} \text{ or } \mathbb{C}, \psi \text{ analytic}$$

with a condition which determines a radius of convergence of the resulting series solution. Based on a geometrical interpretation of Hille's condition, this report generalizes his result to the case when x and y are elements of lattice-normed linear spaces in the sense of Kantorovich. This formulation includes the special cases when x and y belong to Banach or Riesz spaces.

## 0. Introduction

The classical result of Cauchy (1839) states that an implicit function relation in the form

$$y = \psi(x,y) = a_{10}x + \sum_{m+n \geq 2} \sum \frac{a_{mn}}{m!n!} x^m y^n$$

between real or complex variables x,y may, when the series for $\psi(x,y)$ is convergent for sufficiently small x,y, be solved by the method of series substitution in the form

$$y = \varphi(x) = \sum_{n=1}^{\infty} \frac{b_n}{n!} x^n$$

where this series is also convergent for sufficiently small x.

Cauchy's theorem is proved by the method of majorant series (see e.g. Goursat 1904). An associated implicit function equation (comparison equation)

$$Y = \Psi(X,Y) = A_{10}X + \sum_{m+n \geq 2} \sum \frac{A_{mn}}{m!n!} X^m Y^n$$

between real or complex variables X,Y is introduced where $\Psi(X,Y)$ is a majorant to $\psi(x,y)$ in the sense that

$$A_{10} \geq |a_{10}| \; ; \quad A_{mn} \geq |a_{mn}|, \quad m + n \geq 2.$$

It is then easy to show that when the comparison equation is solved in series form as

$$Y = \Phi(X) = \sum_{n=1}^{\infty} \frac{B_n}{n!} X^n$$

that

$$B_n \geq |b_n| \quad n = 0, 1, 2, \ldots \; .$$

Consequently the series for $\varphi(x)$ converges if $|x| \leq X$ and the series for $\Phi(X)$ is convergent. So it is sufficient to consider the solubility by series of the comparison equation which, because of the positivity of $\Psi(X,Y)$, is easier to deal with than the original equation. By choosing a particular majorant function $\Psi(X,Y)$, it is easy to show (Goursat op. cit.) that the corresponding series solution is convergent for sufficiently small values and thus Cauchy's theorem is proved.

Hille (1959), using methods of the theory of functions of a complex variable, strengthened Cauchy's result by showing that if $\Psi(X,Y)$ has a power series which is convergent for all X,Y and, in addition, satisfies some other minor requirements (which are however, apparently not necessary) then the series for $\Phi(X)$ is convergent with radius of convergence $X^*$ given by the solution of the equations

$$Y^* = \Psi(X^*, Y^*), \qquad 1 = \frac{\partial \Psi}{\partial Y}(X^*, Y^*) \qquad \text{(Hille's condition)}.$$

The series for $\varphi(x)$ will then also be convergent when

$$|x| \leq X^*.$$

The significance of Hille's condition is clarified by a geometrical interpretation. Suppose that the graph $\Gamma$ of the comparison equation is constructed for non-negative real X,Y. It is easily shown that, excluding trivial cases, the graph has the form shown below.

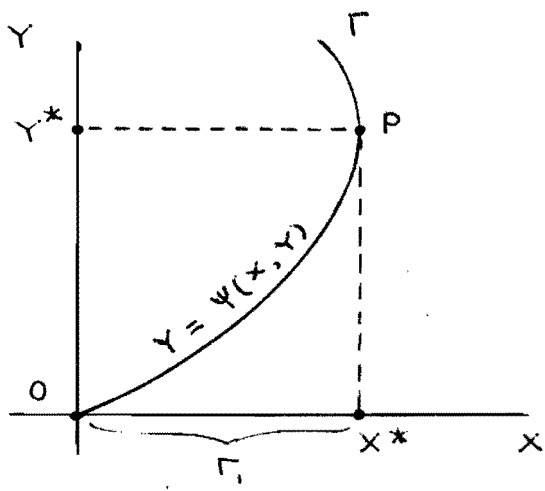

Form of the graph of
the comparison equation.

Starting the origin, the graph curves upward and backward and has a vertical tangent at a point P with coordinates $(X^*, Y^*)$ determined by Hille's condition. The series solution to the comparison equation is thus valid precisely on the arc OP of the graph up to the turning value. The range of convergence of the solution $\Phi(X)$ is, geometrically speaking, the X-projection $\Gamma_1$ of the graph $\Gamma$. The series solution $\varphi(x)$ of the original equation is thus convergent when

$$|x| \leq X, \qquad X \in \Gamma_1.$$

The theorem in this form will be generalized in the present report to the case when x and y are elements of general linear spaces and the variables X and Y belong to associated linear spaces called the norming spaces. These norming spaces are assumed to be Dedekind σ-complete Riesz spaces and an associated completeness condition (called (bk)-completeness) is imposed on the space of elements y. The important special cases when x and y belong to Banach or Riesz spaces are then included in this formulation. The ideas and methods used in making the generalization are due to Kantorovich (1939). The report consists of two parts. The first part develops the theory of analytic functions between linear spaces normed by elements of a Dedekind σ-complete Riesz space. The second part extends the technique of Kantorovich to prove the generalization of Hille's theorem.

# PART I

## PROPERTIES OF MAJORIZED ANALYTIC FUNCTIONS BETWEEN LINEAR SPACES

### 1. Riesz Spaces

The following summarizes the essential points required for the present report. For a complete account see, for example, the books of Vulikh (1967) or Cristescu (1976).

*Definition 1.1.* A *Riesz space* is a real linear space satisfying the following axioms.

(RS1)  If $X_1 \geq X_2$ then $X_1 + X_3 \geq X_2 + X_3$   $\forall X_1, X_2, X_3 \in X$.

(RS2)  If $X_1 \geq X_2$ and $\lambda \geq 0$ then $\lambda X_1 \geq \lambda X_2$   $\forall X_1, X_2 \in X, \lambda \in \mathbb{R}$.

(RS3)  Any $X_1, X_2 \in X$ have a supremum (denoted by $X_1 \vee X_2$).

*Remarks:*

1. It is easily shown that these properties also hold in a dual sense with the relations $\geq, \vee$ on $X$ being replaced by inverse relations $\leq, \wedge$.

2. Elements $X \in X$ satisfying $X \geq 0$ are called *non-negative*. The set of non-negative elements in $X$ will be denoted by $X_+$.

3. Any element $X \in X$ may be written as a difference

$$X = X^+ - X^-$$

of its positive and negative parts $X^+$ and $X^-$ defined by

$$X^+ = X \vee 0, \quad X^- = X \wedge 0$$

From this the absolute value $X$ of any element $X \in X$ may be defined by

$$|X| = X^+ + X^-$$

The absolute value has the usual properties of modulus:

(M1)  $|0| = 0$, $|X| > 0$ if $X \neq 0$   $\forall X \in X$.

(M2) $\quad |X_1 + X_2| \leq |X_1| + |X_2| \qquad \forall\, X_1, X_2 \in X.$

(M3) $\quad |\lambda X| = |\lambda| |X| \qquad \forall\, \lambda \in \mathbb{R},\; X \in X.$

(M4) $\quad ||X_1| - |X_2|| \leq |X_1 - X_2| \qquad \forall\, X_1, X_2 \in X.$

cf. Cristescu p. 74

*Definition 1.2.* A *Dedekind $\sigma$-complete* Riesz space is a Riesz space $X$ satisfying:

(RS4)      An enumerable set $\{X_p\}_{p \in \mathbb{N}}$ of elements of $X$ which has a upper bound has a least upper bound (supremum).

From this follows also the dual property that an enumerable set $\{X_p\}_{p \in \mathbb{N}}$ of elements of $X$ which has a lower bound has a least lower bound (infimum). The supremum and infimum are denoted by the notation $\bigvee_{p \in \mathbb{N}} X_p$ and $\bigwedge_{p \in \mathbb{N}} X_p$ respectively.

A sequence $\{X_p\}_{p \in \mathbb{N}}$ of elements of $X$ is called *monotone increasing* if $X^{(p)} \leq X^{(p+1)}$, $p \in \mathbb{N}$. The property (RS4) enables the convergence of a monotone increasing sequence which is bounded above to be defined, the limit being defined as the supremum $X = \bigvee_{p \in \mathbb{N}} X^{(p)}$. In this case, there is written $X^{(p)} \uparrow X$ as $p \to \infty$. Dually, there is defined convergence of a monotone decreasing sequence which is bounded below.

The convergence of general sequences may be defined by order-convergence (abbreviated (o)-convergence) as in the following definition.

*Definition 1.3.* A sequence $\{X^{(p)}\}_{p \in \mathbb{N}}$ of elements of $X$ is *(o)-convergent* to $X$ as $p \to \infty$ (written $X^{(p)} \overset{(o)}{\to} X$ as $p \to \infty$) if a monotone decreasing sequence of non-negative elements $\{\varepsilon^{(p)}\}_{p \in \mathbb{N}}$ exists such that $\varepsilon^{(p)} \downarrow 0$ as $p \to \infty$ and

$$|X^{(p)} - X| \leq \varepsilon^{(p)} \qquad p \in \mathbb{N},$$

or equivalently,

$$|X^{(p)} - X| \leq \varepsilon^{(n)} \qquad p \geq n \qquad p, n \in \mathbb{N}.$$

Monotone convergence is also order-convergence, e.g. with a monotone increasing sequence $X^{(p)} \uparrow X$ as $p \to \infty$ it is sufficient to take $\varepsilon^{(p)} = X - X^{(p)}$.

The usual properties of limits hold for order-convergence. The important property of the Dedekind σ-complete Riesz space is its completeness with respect to order+convergence as described in the following definition and theorem.

*Definition 1.4.* A sequence $\{X^{(p)}\}_{p \in \mathbb{N}}$ of elements of $X$ is called *(o)-fundamental* if a monotone decreasing sequence of non-negative elements $\{\varepsilon^{(n)}\}_{n \in \mathbb{N}}$ exists such that $\varepsilon^{(n)} \downarrow 0$ as $n \to \infty$ and $|X^{(p)} - X^{(q)}| \leq \varepsilon^{(n)}$, $p,q \geq n$ ($p,q,n \in \mathbb{N}$).

*Theorem 1.1.* A sequence of elements of a Dedekind σ-complete Riesz space is (o)-convergent if and only if it is (o)-fundamental.

Proof: see Cristescu (1976) p. 82.

## 2. Lattice-normed Linear Spaces

Kantorovich (1939) generalized the Banach norm to be an element of a Riesz space. An account of this theory may also be found in the paper of Kantorovich, Vulikh & Pinsker (1951) and in the book of Vulikh (1967). The present section gives a brief account.

*Definition 2.1.* Let $X$ be a linear space and $X'$ a Riesz space. The space $X$ is said to be *lattice-normed* by elements of $X'$ if any $x \in X$ has a norm (modulus, absolute value) denoted by $\|x\|$ such that $\|x\| \in X'$ and the following axioms are satisfied.

(LN1)   $\|x\| \geq 0, \quad \|x\| = 0 \Rightarrow x = 0$   $\forall \, x \in X.$

(LN2)   $\|x_1 + x_2\| \leq \|x_1\| + \|x_2\|,$   $\forall \, x_1, x_2 \in X.$

(LN3)   $\|\lambda x\| = |\lambda| \, \|x\|$   $\forall \, x \in X, \lambda \text{ scalar}.$

If questions of convergence of sequences are involved, it is further necessary to assume that the norming Riesz space is Dedekind $\sigma$-complete. Convergence relative to the norm, called Banach-Kantorovich convergence (abbreviated (bk)-convergence), is then introduced by the definitions:

*Definition 2.2.* A sequence $\{x_p\}_{p \in \mathbb{N}}$ of elements of $X$ is *(bk)-convergent* to $x \in X$ as $p \to \infty$ if $\|x_p - x\| \overset{(o)}{\to} 0$ as $p \to \infty$.

There is written $x_p \overset{(bk)}{\to} x$ as $p \to \infty$ or, alternatively, $x = (bk)\text{-}\lim_{p \to \infty} x_p$. It is easily shown that the (bk)-limit, if it exists, must be unique. The Cauchy property is introduced by:

*Definition 2.3.* A sequence $\{x_p\}_{p \in \mathbb{N}}$ of elements of $X$ is *(bk)-fundamental* if $\exists$ sequence $\{\varepsilon_p\}_{p \in \mathbb{N}} \downarrow 0$ of elements of $X'$ such that

$$\|x_p - x_q\| \leq \varepsilon_n \text{ for all } p, q \geq n.$$

*Definition 2.4.* A linear space, lattice-normed by a Dedekind $\sigma$-complete Riesz space, is said to be *(bk)-complete* if every (bk)-fundamental sequence is (bk)-convergent to some limit.

Following Kantorovich, Vulikh & Pinsker, such a space will be said to be of type $B_K$ (in their paper a further axiom is required for $B_K$ spaces which will not be used here)

*Examples of $B_K$ spaces:*

(a) Let $X$ be a real or complex Banach space and $X' = R$. With $\|x\|$ taken as the Banach norm, (bk)-convergence coincides with ordinary convergence with respect to norm. The space is, by definition, (bk)-complete.

(b) Let $X$ be a Dedekind $\sigma$-complete Riesz space and $X' = X$. With $\|x\|$ taken as absolute value $|x|$ of x in $X$, (bk)-convergence coincides with order-convergence. The space is (bk)-complete (see theorem at end of last section).

(c) Let $X_i$ be a linear space normed by $X_i'$ for $i \in I$. Then the Cartesian product spaces $X = \prod_{i \in I} X_i$ and $X' = \prod_{i \in I} X_i'$ may be formed and $X'$ is by natural definitions, a Riesz space. $x \in X$ is a vector $x = [x_i, i \in I]$ $x_i \in X_i$. $\|x\|$ may be defined as $\|x\| = [\|x_i\|, i \in I]$. It is easy to show that the space is (bk)-complete if the component spaces are.

*Convergent of infinite series:* the usual definitions of convergence of infinite series in Banach or Riesz spaces easily extend to $B_K$ spaces.

Let $\{u_n\}_{n \in \mathbb{N}}$ be a sequence of elements of a $B_K$ space. Introduce the partial sums

$$S_m = \sum_{n=1}^{m} u_n .$$

*Definition 2.5.* The *(bk)-sum* of the infinite series is defined as the limit

$$(bk) - \sum_{n=1}^{\infty} u_n := (bk) - \lim_{m \to \infty} S_m$$

when the limit exists.

This definition applies in particular to the norming Riesz space so that, if $\{U_n\}_{n \in \mathbb{N}}$ is a sequence of elements in the norming space $X'$ then

putting

$$S_m = \sum_{n=1}^{m} U_n$$

there is defined

$$(o) - \sum_{n=1}^{\infty} U_n := (o) - \lim_{m \to \infty} S_m$$

(see Cristescu 1976, p. 83)

Here it is possible to take $U_n = \|u_n\|$ giving:

*Definition 2.6.* An infinite series $\sum u_n$ of elements of a lattice-normed space is said to be *(bk)-absolutely convergent* if $(o) - \sum \|u_n\|$ exists.

An important method of proving convergence of series in a $B_K$ space is by comparison with a series in the norming space as in the following theorem which generalizes the corresponding Riesz space theorem (Cristescu p. 83).

*Theorem 2.1.* Let $\{u_n\}_{n \in \mathbb{N}}$ be a sequence of elements of a $B_K$ space and $\{U_n\}_{n \in \mathbb{N}}$ a sequence of elements in the norming space such that $\|u_n\| \leq U_n$, $n \in \mathbb{N}$. Then if $U = (o) - \sum U_n$ exists, $u = (bk) - \sum u_n$ exists also and $\|u\| \leq U$.

*Proof:* define partial sums $s_m, S_m$ as above. By the Cauchy property, $\exists$ sequence $\{\epsilon_n\} \downarrow 0$ in $X'$ such that

$$\|S_p - S_q\| = \sum_{i=q+1}^{p} U_i \leq \epsilon_n \qquad p \geq q \geq n.$$

Then, for $p \geq q \geq n$,

$$\|s_p - s_q\| = \left\|\sum_{i=q+1}^{p} u_i\right\| \leq \sum_{i=q+1}^{p} \|u_i\| \leq \sum_{i=q+1}^{p} U_i \leq \epsilon_n,$$

and so the result follows from the completeness property of a $B_K$ space.

*Convergence of series of functions:* functions $f: X \to Y$ between lattice--normed linear spaces bring into consideration four related spaces viz. $X$, $Y$ and their norming spaces $X'$, $Y'$. For example, continuity is defined as follows.

*Definition 2.7.* Let $X$, $Y$ be linear spaces lattice-normed by $X'$, $Y'$ respectively. A function $f: X \to Y$ is said to be *(bk)-continuous* if, given $\varepsilon \in X'$, $\varepsilon > 0$, $\exists \eta \in Y'$, $\eta > 0$ such that

$$\|f(x) - f(x')\| < \varepsilon \quad \text{when} \quad \|x - x'\| < \eta \quad \forall x, x' \in X.$$

Consideration of convergence of sequences or series of functions requires the assumption that the space $Y$, into which values are mapped, is (bk)--complete, i.e. $Y$ is a $B_K$ space. Suppose given a sequence $\{u_n(x)\}_{n \in \mathbb{N}}$ from a lattice-normed linear space $X$ to a $B_K$ space $Y$.
The set of values $x \in X$ for which the series is (bk)-convergent will be called the *(bk)-convergence* region $\mathcal{D}$ of the series. There is thus defined for $x \in \mathcal{D}$ the function

$$f(x) = (bk) - \sum_{i=1}^{\infty} u_i(x).$$

Denoting the partial sums by

$$f_n(x) = \sum_{i=1}^{n} u_i(x)$$

there is then a sequence $\{\varepsilon^{(p)}\}_{p \in \mathbb{N}} \downarrow 0$ of elements of $X'$ such that

$$\|f(x) - f_m(x)\| \le \varepsilon^{(n)} \quad \text{if} \quad m \ge n.$$

If the choice of $\varepsilon^{(p)}$ can be made independently of x for $x \in \mathcal{D}_1$ the convergence is said to be *uniform convergence* for $x \in \mathcal{D}_1$. As with the case of ordinary functions of a real variable, there is proved:

*Theorem 1.2.* The (bk)-sum of an infinite series of (bk)-continuous functions which is uniformly convergent in a region is (bk)-continuous in this region.

The theorem above on the use of a comparison series in the norming space is conveniently extended to series of functions in the following way.

*Theorem 1.3.* Let $X$, $Y$ be linear spaces lattice-normed by $X'$, $Y'$ respectively with $Y$ a $B_K$ space. Let $\{u_n(x)\}_{n \in \mathbb{N}}$ and $\{U_n(X)\}_{n \in \mathbb{N}}$ be sequences of functions $X \to Y$ and $X' \to Y'$ respectively, satisfying

$$\|u_n(x)\| \leq U_n(x) \quad \text{when} \quad \|x\| \leq X, \quad n \in \mathbb{N}.$$

Then if $\sum U_n(X)$ has the (o)-convergence set D', the series $\sum u_n(x)$ is (bk)-convergent for $\|x\| \leq X \in D'$.

Here the functions $U_n(X)$ will often be *isotone* for $X \geq 0$, i.e.

$$U_n(X) \leq U_n(X') \quad \text{when} \quad 0 \leq X \leq X'$$

Define an *order-interval* ((o)-interval) $<X_1, X_2>$ for any $X_1, X_2 \in X'$ as the set

$$<X_1, X_2> = \{X \in X' \mid X_1 \leq X \leq X_2\}$$

and an *order-star* (abbreviated (o)-star) as a set $S \subset X'$ which contains, with any X, the whole (o)-interval $<0, X>$. Then it is clear that:

*Corollary.* If, in the last theorem, $U_n(X)$ $n \in \mathbb{N}$ is isotone for $X \geq 0$ then the (o)-convergence region D' of $\sum U_n(X)$ is an (o)-star and $\sum u_n(x)$ will (bk)-converge when $\|x\| \in D'$.

## 3. Majorized Analytic Functions

This section develops a theory of analytic functions in lattice-normed linear spaces. There does not appear to be any account of this theory in the literature although Kantorovich, Vulikh & Pinsker (1951) had observed that the construction of such a theory is possible.*

*Multilinear functions*: let $X$ and $Y$ be real or complex linear spaces. A *multilinear function* from $X$ to $Y$ of degree n (*n-linear function*) is a function $a_n(x_1,...,x_n)$ of $x_1,...,x_n \in X$ with values in $Y$ which is linear in each of $x_1,...,x_n$. The set of multilinear functions of a given degree forms a linear space under addition and scalar multiplication which will be denoted by $L^n(X,Y)$. A multilinear function of degree 0 is defined as a constant value in $Y$ so that $L^0(X,Y)$ is $Y$ itself.

Of special importance are completely symmetrical multilinear functions.

*Definition 3.1.* $a_n \in L^n(X,Y)$ is said to be *completely symmetrical* if

$$a_n(x_{i_1},...,x_{i_n}) = a(x_1,...,x_n)$$

for all permutations $i_1,...,i_n$ of $1,...,n$. From any $a_n \in L^n(X,Y)$, a completely symmetrical element $a'_n \in L^n(X,Y)$ can be constructed as

$$a'_n(x_1,...,x_n) = \frac{1}{n!} \sum_{\text{perms}} a_n(x_{i_1},...,x_{i_n})$$

where the sum on the right is over all permutations $i_1,...,i_n$ of $1,...n$. $a'_n$ will be referred to as the *symmetrised form* of $a_n$.

*Polynomials*: a function $u_n : X \to Y$ may be constructed from any $a_n \in L^n(X,Y)$ by

$$u_n(x) = a_n(x,...,x) \qquad x \in X.$$

It is seen to be homogeneous of degree n in the sense that

$$u_n(\lambda x) = \lambda^n u_n(x) \qquad x \in X, \lambda \text{ scalar}$$

and may be called a *monomial* of degree n in x.

---

* See p. 106 of the Amer. Math. Soc. Trans.

If $a_n' \in L^n(X,Y)$ is the symmetrized form of $a_n \in L^n(X,Y)$ then clearly

$$a_n'(x,\ldots,x) = a_n(x,\ldots,x) \qquad \forall x \in X,$$

i.e. $a_n'$ and $a_n$ correspond to the same monomial. It may be shown (Mazur & Orlicz 1935, Hille & Phillips 1948) that, with certain additional assumptions, the symmetrized form of a multilinear function may be reconstructed from its monomial. Thus there is a one-to-one correspondence between monomials and completely symmetrical multilinear functions.

A *polynomial of degree n* may now be defined as a sum of monomials of degrees $0, 1, \ldots, n$. In forming polynomials, the completely symmetrical forms of the corresponding multilinear functions are normally used.

*Power series*: given a sequence $\{a_n\}_{n \in \mathbb{N}}$, $a_n \in L^n(X,Y)$, a *power series* may be defined formally as a sum

$$\sum_{n=1}^{\infty} \frac{1}{n!} a_n(x,\ldots,x)$$

where the $n!$ factor has been put in for convenience. The functions $a_n$ are called the *coefficients* of the power series. Summation may also start at $n = 0$ if a constant term $a_0$ is added.

The (bk)-convergence of such a power series is defined using the theory of the last section. Denoting the region of (bk)-convergence by $\mathcal{D}$ there is defined for $x \in \mathcal{D}$ a function

$$y = f(x) = (bk) - \sum_{n=1}^{\infty} \frac{1}{n!} a_n(x,\ldots,x)$$

which will be called a *(bk)-analytic function* $f : X \to Y$.

*Power series in the norming space*: the previous definitions apply, in particular, also to the norming spaces $X', Y'$ of $X, Y$. In this way are defined (o)-analytic functions

$$Y = F(X) = (o) - \sum_{n=1}^{\infty} \frac{1}{n!} A_n(X,\ldots,X)$$

within a convergence region $\mathcal{D}' \subset X'$.

*Definition 3.2.* A power series in the norming space defined by coefficients $A_n$, $n \in \mathbb{N}$ is said to be of *positive type* if for $n \in \mathbb{N}$

$$A_n(X_1,\ldots,X_n) \geq 0 \quad \text{when} \quad X_1,\ldots,X_n \geq 0 \quad \forall X_1,\ldots,X_n \in X'.$$

Power series of positive type are of particular importance because their convergence properties are simple. The partial sums of such a power series are positive and monotone increasing. The region of convergence is consequently the set on which these partial sums are bounded above. Since the individual terms of a power series of positive type are clearly isotone functions of X when $X \geq 0$ there follows (cf. last section)

*Lemma 3.1.* The region of convergence of a power series of positive type is an (o)-star.

*Majorized power series:* consider now the relation between power series between general lattice-normed linear spaces and power series in the norming spaces.

*Definition 3.3.* A power series $f : X \to Y$ defined by coefficients $a_n$, $n \in \mathbb{N}$ is said to be *majorized* by a power series of positive type $F : X' \to Y'$ defined by coefficients $A_n$, $n \in \mathbb{N}$ if for $n \in \mathbb{N}$

$$\| a_n(x_1,\ldots,x_n) \| \leq A_n(X_1,\ldots,X_n)$$

for all $x_1,\ldots,x_n \in X$ such that $\|x_i\| \leq X_i$, $i = 1,\ldots,n$.

The majorant defined in this way is a direct generalization of the original Cauchy majorant for ordinary power series. It may thus be also called the *Cauchy majorant.* If f is majorized by F it is convenient to write $f \ll F$ in a notation due to Poincaré. In this case $a_n \ll A_n$ $n \in \mathbb{N}$ also.

Note that if the coefficients $a_n$, $n \in \mathbb{N}$ are majorized by $A_n$, $n \in \mathbb{N}$ then the same will be true for the completely symmetrized forms of these coefficients.

The following theorem results directly from the discussion of convergence given in the last section.

*Theorem 3.1.* Suppose that, as in the last definition, $f \ll F$ where $F$ has region of convergence $\mathcal{D}'$. Then

(a) the series for $f$ is (bk)-absolutely convergent for $x \in \mathcal{D} = \{x \in X \mid \|x\| \in \mathcal{D}'\}$
(b) the series for $f$ is uniformly (bk)-convergent for $x \in \mathcal{D}_X = \{x \in X \mid \|x\| \leq X \in \mathcal{D}'\}$
(c) when $x \in \mathcal{D}$ the inequality $\|f(x)\| \leq F(\|x\|)$ holds.

In the case when the coefficients $a_n$, $n \in \mathbb{N}$, of a power series $X \to Y$ are (bk)-continuous functions in each of their variables, the monomials $a_n(x,\ldots,x)$ formed from them are also (bk)-continuous functions $X \to Y$. Hence, using the uniform convergence property proved in the last theorem there follows:

*Corollary.* If in the last theorem, the coefficients $a_n, A_n$ are (bk)-continuous ((o)-continuous) multilinear functions, then $f(x), F(X)$ are (bk)-continuous((o)-continuous) for $\|x\| \leq X \in \mathcal{D}'$.

*Regular power series in Riesz space:* the case when $X$ and $Y$ are Riesz spaces lattice-normed by themselves is of special importance. Assume $X = X'$, $Y = Y'$, $\|x\| = |x|$, $\|y\| = |y|$.

*Definition 3.4.* $a_n \in L^n(X,Y)$ is called *positive* (written $a_n \geq 0$) if

$$a_n(x_1,\ldots,x_n) \geq 0 \quad \text{when} \quad x_1,\ldots,x_n \geq 0 \qquad \forall x_1,\ldots,x_n \in X.$$

The set of positive elements of $L^n(X,Y)$ will be denoted by $L^n_+(X,Y)$.

*Lemma 3.2.* If $a_n \in L^n_+(X,Y)$ then

$$|a_n(x_1,\ldots,x_n)| \leq a_n(|x_1|,\ldots,|x_n|) \qquad \forall x_1,\ldots,x_n \in X.$$

*Proof.* see Cristescu 1971, p. 201.

The definition of positivity gives a partial ordering of the linear space $L^n(X,Y)$ by defining $b_n \geq a_n$ to mean $b_n - a_n \geq 0$. (see Cristescu loc. cit.)

*Definition 3.5.* A multilinear function $a_n \in L^n(X,Y)$ is called *regular* if it can be represented as

$$a_n = a_n^+ - a_n^-, \quad a_n^+, a_n^- \in L_+^n(X,Y).$$

It is then possible to define an absolute value by

$$|a_n| = a_n^+ + a_n^-.$$

(Cristescu, loc. cit.)

*Lemma 3.3.* If $a_n \in L^n(X,Y)$ is regular then

$$|a_n(x_1,\ldots,x_n)| \leq |a_n|(|x_1|,\ldots,|x_n|) \quad \forall x_1,\ldots,x_n \in X.$$

*Proof.*

$$|a_n(x_1,\ldots,x_n)| \leq |a_n^+(x_1,\ldots,x_n) - a_n^-(x_1,\ldots,x_n)|$$

$$\leq |a_n^+(x_1,\ldots,x_n)| + |a_n^-(x_1,\ldots,x_n)|$$

$$\leq a_n^+(|x_1|,\ldots,|x_n|) + a_n^-(|x_1|,\ldots,|x_n|)$$

$$= |a_n|(|x_1|,\ldots,|x_n|).$$

*Definition 3.6.* A power series from a Riesz space to a Dedekind $\sigma$-complete Riesz space will be called *regular* if all its coefficients are regular.

A majorant always exists for a regular power series, namely the one with coefficients $A_n = |a_n|$, $n \in \mathbb{N}$. Any power series with coefficients $A_n \in L_+^n(X,Y)$ such that $A_n \geq |a_n|$ will then also be a majorant. Of these majorants, that with coefficients $|a_n|$ will have greatest region of convergence.

*An inequality satisfied by the majorant:* returning to the general case of lattice-normed linear spaces, this section will end by proving an inequality due to Kantorovich (1939) which will be used in the second part of this report.

*Lemma 3.4.* In the previous notation, if $f \ll F$ and $F$ has region of convergence $\mathcal{D}'$ then

$$\| f(x + h) - f(x) \| \leq F(X + H) - F(X)$$

when

$$\| x \| \leq X, \quad \| h \| \leq H, \quad X, X + H \in \mathcal{D}' .$$

*Proof.* The stated inequality is proved first for a coefficient pair $a_n, A_n$ where $a_n \ll A_n$.

$$a_n(x + h, \ldots, x + h) - a_n(x, \ldots, x)$$

$$= \sum_{i=1}^{n} \left\{ \begin{array}{l} a_n(\underbrace{x + h, \ldots, x + h}_{(n-i+1)}, \underbrace{x, \ldots, x}_{(i-1)}) \\ -a_n(\underbrace{x + h, \ldots, x + h}_{(n-i)}, \underbrace{x, \ldots, x}_{i}) \end{array} \right\}$$

$$= \sum_{i=1}^{n} a_n(\underbrace{x + h, \ldots, x + h}_{(n-i)}, h, \underbrace{x, \ldots, x}_{(i-1)}) .$$

On taking norms there follows

$$\| a_n(x + h, \ldots, a + h) - a_n(x, \ldots, x) \|$$

$$\leq \sum_{i=1}^{n} A_n(\underbrace{X + H, \ldots, X + H}_{(n-i)}, H, \underbrace{X, \ldots, X}_{(i-1)})$$

$$= A_n(X + H, \ldots, X + H) - A_n(X, \ldots, X)$$

on reversing the previous algebra. This proves the stated inequality for a pair of coefficients. Now follows

$$\|f(x + h) - f(x)\|$$

$$\leq \sum_{n=1}^{\infty} \frac{1}{n!} \|a_n(x + h, \ldots, x + h) - a_n(x, \ldots, x)\|$$

$$\leq \sum_{n=1}^{\infty} \frac{1}{n!} (A_n(X + H, \ldots, X + H) - A_n(X, \ldots, X))$$

$$= F(X + H) - F(X)$$

as required. The operations on the right here are justified by the absolute convergence of the series.

## 4. Use of the Cauchy Majorant in Series Substitution

*Direct series substitution:* suppose that $X$, $Y$, $Z$, are three lattice-normed linear spaces and $f : X \to Y$ and $g : Y \to Z$ are (formal) power series given by

$$y = f(x) = \sum_{n=0}^{\infty} \frac{1}{n!} a_n(x,\ldots,x),$$

$$z = g(y) = \sum_{n=0}^{\infty} \frac{1}{n!} b_n(y,\ldots,y).$$

Substitution for y gives a series for $h = g \circ f : X \to Z$ :

$$z = h(x) = \sum_{n=0}^{\infty} \frac{1}{n!} c_n(x,\ldots,x).$$

It is easily verified that, for completely symmetrical coefficients, the following relations hold:

$$c_1(x) = b_1 \circ a_1(x)$$

$$c_2(x_1, x_2) = b_1 \circ a_2(x_1, x_2) + b_2(a_1(x_1), a_2(x_2))$$

$$\ldots \quad \ldots \quad \ldots \quad \ldots \quad \ldots \quad \ldots \quad \ldots$$

$$c_n(x_1, \ldots, x_n) =$$

$$\sum_{r=1}^{n} \sum_{\substack{n_1 + \ldots + n_r \\ = n}} \sum_{\pi(n_1, \ldots, n_r)} b_r(a_{n_1}(x_{i_1}, \ldots, x_{i_{n_1}}), \ldots, a_{n_r}(x_{i_{n-n_r+1}}, \ldots, x_{i_n}))$$

where the last sum on the right is over all partitions $\pi(n_1, \ldots, n_r)$ of the set $\{1, \ldots, n\}$ of suffices into r non-null disjoint sets with $n_1, \ldots, n_r$ elements respectively of the type

$$\{1, \ldots, n\} = \{i_1, \ldots, i_{n_1}\} \cup \ldots \cup \{i_{n-n_r+1}, \ldots, i_n\}.$$

The formulae for $c_n$ given here is a generalization of the well-known di Bruno formula which applies to substitution of ordinary power series (Riordan 1968).

Suppose now that $f \ll F$ and $g \ll G$ where $F : X' \to Y'$ and $G : Y' \to Z'$ are given by

$$Y = F(X) = \sum_{n=1}^{\infty} \frac{1}{n!} A_n(X,\ldots,X),$$

$$Z = G(Y) = \sum_{n=1}^{\infty} \frac{1}{n!} B_n(Y,\ldots,Y),$$

where the prefix (o) has been omitted since only formal algebraic properties are in question at the moment. Substitution gives the series for $H = G \circ F : X' \to Z'$:

$$Z = H(X) = \sum_{n=1}^{\infty} \frac{1}{n!} C_n(X,\ldots,X)$$

where the coefficients $C_n$ are related to the coefficients $A_n$ and $B_n$ by formulae similar to those above.

The basic property of the Cauchy majorant is expressed by:

*Lemma 4.1.* If $f \ll F$ and $g \ll G$ then $g \circ f \ll G \circ F$.

*Proof.* It is required to show that for $n \in \mathbb{N}$,

$$\| c_n(x_1,\ldots,x_n) \| \leq C_n(X_1,\ldots,X_n) \quad \text{when} \quad \|x_i\| \leq X_i \quad i = 1,\ldots,n .$$

This inequality may be proved by induction using the above algebraic formulae. For $n = 1$:

$$\| c_1(x) \| = \| b_1(a_1(x)) \| \leq B_1(\| a_1(x)\|)$$

$$\leq B_1(A_1(X)) = C_1(X) .$$

Assuming that the inequality is valid for the values $0, 1, \ldots, n-1$, there follows:

$$\|c_n(x_1,\ldots,x_n)\| \le$$

$$\le \sum_{r=1}^{n} \sum_{\substack{n_1+\ldots+n_r \\ =n}} \sum_{\pi(n_1,\ldots,n_r)} \|b_r(a_{n_1}(x_{i_1},\ldots,x_{i_{n_1}}),\ldots)\|$$

$$\le \sum_{r=1}^{n} \sum_{\substack{n_1+\ldots+n_r \\ =n}} \sum_{\pi(n_1,\ldots,n_r)} B_r(\|a_{n_1}(x_{i_1},\ldots,x_{i_{n_1}})\|,\ldots)$$

$$\le \sum_{r=1}^{n} \sum_{\substack{n_1+\ldots+n_r \\ =n}} \sum_{\pi(n_1,\ldots,n_r)} B_r(A_{n_1}(X_{i_1},\ldots,X_{i_{n_1}}),\ldots)$$

$$= C_n(X_1,\ldots,X_n)$$

completing the induction.

Using this lemma and previous theorems on majorant series there results the theorem on series substitution:

*Theorem 4.1.* Let $f : X \to Y$ and $g : Y \to Z$ be power series between lattice-normed linear spaces $X$, $Y$, $Z$ where $Y$ and $Z$ are $B_K$ spaces. Suppose f and g are majorized by $F : X' \to Y'$ and $G : Y' \to Z'$ respectively, with regions of convergence $\mathcal{D}'$, $\mathcal{D}''$. Putting $\mathcal{D}'_1 = \{X \in \mathcal{D}' | F(X) \in \mathcal{D}''\}$ there follows

(a) the series for $G \circ F$ is (o)-convergent for $X \in \mathcal{D}'_1$
(b) the series for $g \circ g$ is (bk)-convergent for $\|x\| \in \mathcal{D}'_1$.

*Series substitution for analytic implicit function equations:* now consider the solution by series of an implicit function equation

$$y = \psi(x,y)$$

where $\psi : X \times Y \to Y$ is represented by the (formal) series

$$\psi(x,y) = \sum_{m+n \geq 2} \frac{1}{m!n!} a_{mn}(\underbrace{x,\ldots,x}_{m}; \underbrace{y,\ldots,y}_{n}).$$

Here the coefficients $a_{mn}$ may be assumed to be multilinear functions of $m + n$ variables which are completely symmetric in the first m and last n variables.

It will be assumed that the equation has already be solved with respect to the linear term in y so that the linear term in y on the right is absent, i.e. $a_{01} = 0$. The equation thus has the form

$$y = a_{10}x + \psi_1(x,y)$$

where $\psi_1$ contains only nonlinear terms, i.e. terms of $\psi$ of degree $m + n \geq 2$.

Solution of the equation by series substitution means determining a series $\varphi : X \to Y$:

$$y = \varphi(x) = \sum_{n=0}^{\infty} \frac{1}{n!} b_n(x,\ldots,x)$$

which formally satisfies the equation. On substitution the following relations are found for the completely symmetrical coefficients $b_n$:

$$b_1(x) = a_{10}(x)$$

$$b_2(x_1,x_2) = a_{20}(x_1,x_2)$$
$$+ a_{11}(x_1, b_1(x_2)) + a_{11}(b(x_1), x_2)$$
$$+ a_{02}(b_1(x_1), b_1(x_2))$$

... ... ... ... ...

$$b_n(x_1,\ldots,x_n) =$$

$$\sum_{r+m=n} \sum_{\substack{n_1+\ldots+n_r \\ =m}} \sum_{\pi(m,n_1,\ldots,n_r)} a_{mr}(x_{i_1},\ldots,x_{i_m}; b_{n_1}(x_{i_{m+1}},\ldots,x_{m+n_1}),\ldots,b_{n_r}(x_{i_{n-n_r+1}},\ldots,x_{i_n}))$$

where the last summation on the right is over all partitions $\pi(m,n_1,\ldots,n_r)$ of the set of suffices $1,\ldots,n$ into $1 + r$ non-null disjoint sets with $m,n_1,\ldots,n_r$ elements respectively of the type

$$\{1,\ldots,n\} = \{i_1,\ldots,i_m\} \cup [\{i_{m+1},\ldots,i_{m+n_1}\} \cup \ldots$$
$$\ldots \cup \{i_{n-n_r+1},\ldots,i_n\}].$$

The above equations determine the coefficients recursively as completely symmetrical multilinear functions. A formal solution of the implicit equation is thus uniquely determined.

Now consider the equation

$$Y = \Psi(X,Y)$$

in the norming spaces where $\Psi : X' \times Y' \to Y'$ has a series representation

$$\Psi(X,Y) = \sum_{m+n \geq 1} \frac{1}{m!n!} A_{mn} (\underbrace{X,\ldots,X}_{m}; Y,\ldots,Y)$$

with $A_{01} = 0$. This equation will be called the *comparison equation* to the previous implicit function equation. It may be solved in a similar way as a series

$$Y = \Phi(X) = \sum_{n=1}^{\infty} \frac{1}{n!} B_n(X,\ldots,X)$$

where the $B_n$ are determined recursively by equations similar to those for the $b_n$. The idea of Cauchy is now expressed by the following lemma:

*Lemma 4.2.* If $\psi \ll \Psi$ then $\varphi \ll \Phi$.

*Proof.* It is required to show that for $n \in \mathbb{N}$,

$$\| b_n(x_1,\ldots,x_n) \| \leq B_n(X_1,\ldots,X_n) \quad \text{when} \quad \|x_i\| \leq X_i \quad i = 1,\ldots,n.$$

This inequality is proved by induction on $n$ in a similar way to the previous lemma on series substitution. The proof may be omitted.

*Convergence*: with the same notation, suppose:

(a)  the series for $\Psi$ is (o)-convergent for $(X,Y) \in \mathcal{E}'$ where $\mathcal{E}'$ is an (o)-star in $X'_+ \times Y'_+$.

(b)  the series for $\Phi$ is (o)-convergent for $X \in \mathcal{D}'$ where $\mathcal{D}'$ is an (o)-star in $X'_+$.

The following lemma is a straightforward extension of theorem 4.1 on substitution of power series.

*Lemma 4.3.* Suppose that the above convergence conditions are satisfied. Then

(a)  series substitution of the power series of $\Phi$ into the power series for $\Psi$ is permissible and leads to an (o)-convergent power series for $\Psi(X,\Phi(X))$ when

$$X \in \mathcal{D}' \cap \{X \in X'_+ | \exists \, Y \in Y'_+ : (X,Y) \in \mathcal{E}' \,\&\, 0 \leq \Phi(X) \leq Y\}$$

If also $\Psi \gg \psi$ then $\Phi \gg \varphi$ by lemma 4.2 and there follows

(b)  if, for such a pair $X,Y$, $x \in X$, $y \in Y$ satisfy $\|x\| \leq X$, $\|y\| \leq Y$, then series substitution of the power series for $\varphi$ into the power series for $\psi$ is permissible and leads to a (bk)-convergent power series for $\psi(x,\varphi(x))$.

The final conclusion of this section may be stated in the form of a theorem.

*Theorem 4.2.* Consider the analytic equations $y = \psi(x,y)$, $Y = \Psi(X,Y)$ where $\psi \ll \Psi$, it being assumed that $\psi, \Psi$ have no linear term in $y, Y$ respectively. Let $y = \varphi(x)$, $Y = \Phi(X)$ be the series solutions of these equations.
If $\Psi$ and $\Phi$ have convergence as given in conditions (a), (b) above then

(a)  $Y = \Phi(X)$ is an (o)-analytic solution of $Y = \Psi(X,Y)$ with region of convergence $X \in \mathcal{D}' \cap \{X | (X,\Phi(X)) \in \mathcal{E}'\}$.

(b)  $y = \varphi(x)$ is a (bk)-analytic solution of $y = \psi(x,y)$ with region of convergence $x \in \{x | \|x\| \in \mathcal{D}' \,\&\, (\|x\|, \Phi(\|x\|)) \in \mathcal{E}'\}$.

In this theorem it is necessary to assume that $\Phi(X)$ has a convergent power series. The next part of this report is devoted to proving that this is so.

PART II

THE APPLICATION OF KANTOROVICH'S THEORY TO IMPLICIT FUNCTION EQUATIONS

## 5. Kantorovich's Method of Successive Approximation

In his 1939 paper, Kantorovich used the lattice-norm in connection with the convergence of the method of successive approximation for equations of the type

$$x = f(x) \qquad x \in X.$$

For this purpose he introduced an auxiliary equation of similar form

$$X = F(X) \qquad X \in X'$$

defined in the norming space (the *comparison equation*). The norming space was assumed to be of type $B_K$ and the function $F(X)$ was assumed to satisfy the following conditions for values of $X$ on an order interval $\mathcal{D}' = \langle 0, \tilde{X} \rangle$.

K1  $F(X)$ is defined for $X \in \mathcal{D}'$.

K2  If $F$ is continuous with respect to monotone increasing sequences on $\mathcal{D}'$.

K3  $F$ is isotone on $\mathcal{D}'$ i.e. if $X, X' \in \mathcal{D}'$, $X \leq X'$ then $F(X) \leq F(X')$.

K4  $F(0) \geq 0$.

K5  $\tilde{X} \geq F(\tilde{X})$.

The relation between $F$ and $f$ was assumed to be governed by :

M1  $f(x)$ is defined for $x \in \mathcal{D} = \{x \in X \mid \|x\| \in \mathcal{D}'\}$.

M2  $\|f(0)\| \leq F(0)$.

M3  $\|f(x + h) - f(x)\| \leq F(X + H) - F(X)$ when

$\|x\| \leq X$, $\|h\| \leq H$, and $X, X + H \in \mathcal{D}'$.

Note that these conditions are satisfied by a (bk)-analytic function and its majorant on any order interval of convergence, provided that the coefficients of the majorant are continuous.

From K3 and K4 there follows immediately

K6  $F(X) \geq 0$  when  $X \geq 0$, $X \in \mathcal{D}'$.

From M3 follows, on changing x, X, h, H to 0, 0, x, X respectively.

$$\| f(x) - f(0) \| \leq F(X) - F(0) \quad \text{when } \|x\| \leq X \in \mathcal{D}'$$

from which follows, using M2,

$$\| f(x) \| \leq \| f(0) \| + \| f(x) - f(0) \| \leq F(0) + (F(X) - F(0)) = F(X)$$

so that there is satisfied

M4  If $\|x\| \leq X \in \mathcal{D}'$ then $\| f(x) \| \leq F(X)$.

*Lemma 5.1.* If conditions K1 - K5 and M1 - M3 are satisfied, the functions f, F, map the sets $\mathcal{D}$, $\mathcal{D}'$ respectively, into themselves.

*Proof.* If $X \in \mathcal{D}'$ then by K6, K3, K5,

$$0 \leq F(X) \leq F(\widetilde{X}) \leq \widetilde{X}$$

if $x \in \mathcal{D}$, $\|x\| \leq \widetilde{X}$ so by M4, K5

$$\| f(x) \| \leq F(\widetilde{X}) \leq \widetilde{X} .$$

From this lemma it follows that the iterative process

$$x^{(0)} = 0$$
$$x^{(p+1)} = f(x^{(p)}) , \quad p \in \mathbb{N}$$

and the associated process

$$X^{(0)} = 0$$
$$X^{(p+1)} = F(X^{(p)}) , \quad p \in \mathbb{N}$$

are well-defined and give sequences $\{x^{(p)}\}_{p \in \mathbb{N}}$, $\{X^{(p)}\}_{p \in \mathbb{N}}$ lying in $\mathcal{D}$, $\mathcal{D}'$ respectively.

*Lemma 5.2.* The sequence $\{X^{(p)}\}_{p \in \mathbb{N}}$ is monotone increasing on $\mathcal{D}'$.

*Proof.* Induction. For $p = 0$,

$$X^{(1)} = F(0) \geq 0 = X^{(0)} .$$

For p > 0, assuming that

$$X^{(p)} \geq X^{(p-1)} \geq 0$$

there follows

$$X^{(p+1)} - X^{(p)} = F(X^{(p)}) - F(X^{(p-1)}) \geq 0$$

completing the induction.

*Lemma 5.3.* The following relations hold between the iterative process on $\mathcal{D}$, $\mathcal{D}'$:

(a) $\quad \|x^{(p)}\| \leq X^{(p)}, \quad p \in \mathbb{N}$.

(b) $\quad \|x^{(p+1)} - x^{(p)}\| \leq X^{(p+1)} - X^{(p)}, \quad p \in \mathbb{N}$.

*Proof.* Induction:

(a), p = 0:

$$\|x^{(0)}\| = 0 = X^{(0)}.$$

For p > 0, assuming that

$$\|x^{(p)}\| \leq X^{(p)}$$

there follows

$$\|x^{(p+1)}\| \leq \|f(x^{(p)})\| \leq F(\|x^{(p)}\|) \leq F(X^{(p)}) = X^{(p+1)}$$

proving (a).

(b), p = 0: since $x^{(0)} = 0$, $X^{(0)} = 0$,

$$\|x^{(1)} - x^{(0)}\| = \|x^{(1)}\| = \|f(x^{(0)})\|$$

$$\leq F(X^{(0)}) = X^{(1)} = X^{(1)} - X^{(0)}.$$

For p > 0, assuming that

$$\|x^{(p)} - x^{(p-1)}\| \leq X^{(p)} - X^{(p-1)},$$

there follows

$$\|x^{(p+1)} - x^{(p)}\| = \|f(x^{(p)}) - f(x^{(p-1)})\|$$
$$\leq F(X^{(p)}) - F(X^{(p-1)})$$
$$= X^{(p+1)} - X^{(p)}.$$

The main result of Kantorovich may now be stated.

*Theorem 5.1.* (Kantorovich) The iterates $x^{(p)}$, $X^{(p)}$, $p \in \mathbb{N}$ satisfy

(a)   $X^{(p)} \uparrow X$ where $X \in \mathcal{D}'$ and $X$ satisfies $X = F(X)$ .
(b)   $x^{(p)} \overset{(bk)}{\to} x$ as $p \to \infty$ where $x \in \mathcal{D}$, $\|x\| \leq X$ and $x = f(x)$.

*Proof.*

(a) From lemma 5. $X^{(p)} \uparrow X$ where $0 \leq X \leq \widetilde{X}$. Since $F$ is assumed continuous for increasing sequences, $X = F(X)$.

(b) If $p \geq q \geq n$, $p,q,n \in \mathbb{N}$,

$$\|x^{(p)} - x^{(q)}\| \leq \|x^{(p)} - x^{(p-1)}\| + \ldots + \|x^{(q+1)} - x^{(q)}\|$$
$$\leq (X^{(p)} - X^{(p-1)}) + \ldots + (X^{(q+1)} - X^{(q)})$$
$$\leq X^{(p)} - X^{(q)} \overset{(o)}{\to} 0 \quad \text{as } n \to \infty$$

Consequently $\{x^{(p)}\}_{p \in \mathbb{N}}$ is a fundamental sequence and so (bk)-convergent to a limit $x \in X$. Since, by lemma 5.3, $\|x^{(p)}\| \leq \|X^{(p)}\|$ it follows that $\|x\| \leq X$ so that $x \in \mathcal{D}$.

Now from

$$x - x^{(q)} = \sum_{p=q+1}^{\infty} (x^{(p)} - x^{(p-1)})$$

it follows that

$$\|x - x^{(q)}\| \leq \sum_{p=q+1}^{\infty} \|x^{(p)} - x^{(p-1)}\|$$
$$\leq \sum_{p=q+1}^{\infty} (X^{(p)} - X^{(p-1)}) = X - X^{(q)} .$$

Now substituting $x^{(p)}$, $x - x^{(p)}$ for $x$, $h$, in condition M3 it follows that for $X^{(p)}$, $X \in \mathcal{D}$,

$$\| f(x) - f(x^{(p)}) \| \leq F(X) - F(X^{(p)})$$

Since $F(X^{(p)}) \uparrow F(X)$ there follows $f(x^{(p)}) \overset{(bk)}{\to} f(x)$ as $p \to \infty$.

As a corollary to the last theorem it is seen that in Kantorovich's conditions K1 - K5, the order interval $<0,\tilde{X}>$ may be replaced by $<0,X>$, the iterative sequence $\{X^{(p)}\}_{p \in \mathbb{N}}$ lying entirely on the interval $<0,X>$.

A more convenient statement of the Kantorovich conditions for the present purpose comes about by taking the region $\mathcal{D}'$ of definition of F to be an arbitrary order-star in $X'_+$. $\mathcal{D}$ may be defined from $\mathcal{D}'$ by M1. Condition K5 is changed to

K5 $\exists \tilde{X} \in \mathcal{D}' : \tilde{X} \geq F(\tilde{X})$ .

Under this slight modification, $<0,\tilde{X}>$ lies in $\mathcal{D}'$ so that the previous conclusions still hold. Taking $\mathcal{D}'$ as an order star is more natural when considering analytic functions and it also permits the statement of the following minimal property:

*Theorem 5.2.* Suppose that Kantorovich's conditions hold in the form K1-K5, M1-M3 when $\mathcal{D}'$ is an order star in $X'_+$. Suppose that a solution $X' \in \mathcal{D}'$ of the equation $X = F(X)$ exists. Then the solution X corresponding to the previous iterative process also exists and $X' \geq X$.

*Proof.* The order interval $<0,X'>$ lies in $\mathcal{D}'$ and so may be used instead of the order interval $<0,\tilde{X}>$ in the previous results. Hence the iterative process lies in $<0,X'>$ and tends to a limit X satisfying $X \leq X'$.

## 6. The Implicit Function Theorem under Kantorovich Conditions

Let $X$ and $Y$ be linear spaces lattice-normed by $X'$ and $Y'$ respectively. Assume that $Y$ is a $B_k$ space.

Consider the implicit function equation

$$y = \psi(x,y) \qquad x \in X, \; y \in Y,$$

and the comparison equation

$$Y = \Psi(X,Y) \qquad X \in X'_+, \; Y \in Y'_+ .$$

Suppose that the following modified Kantorovich conditions are satisfied:

K1' $\Psi(X,Y)$ is defined when $(X,Y) \in \mathcal{E}'$ where $\mathcal{E}'$ is an order star in in $X_+ \times Y_+$ ,

K2 $\Psi(X,Y)$ is continuous with respect to monotone increasing sequences in $\mathcal{E}'$,

K3' $\Psi$ is isotone on $\mathcal{E}'$ with respect to $(X,Y)$.

K4' $\Psi(0,0) \geq 0$ ,

K5' For a non-null set of values of $X$ in $X'_+$, $\exists \, Y$ such that $(X,Y) \in \mathcal{E}'$ and $Y \geq \Psi(X,Y)$,

and that $\Psi(X,Y)$ and $\psi(x,y)$ are related by the conditions

M1' $\psi(x,y)$ is defined for $x \in \mathcal{D} = \{(x,y) \in X \times Y \mid (\|x\|, \|y\|) \in \mathcal{E}'\}$,

M2' $\|\psi(0,0)\| \leq \Psi(0,0)$,

M3' $\|\psi(x,y+k) - \psi(x,y)\| \leq \Psi(X,Y+K) - \Psi(X,Y)$
when $\|x\| \leq X$, $\|y\| \leq Y$, $\|k\| \leq K$; $(X,Y), (X,Y+K) \in \mathcal{E}'$.

From these conditions the two further conditions

K6' $\Psi(X,Y) \geq 0$ when $X \geq 0$, $Y \geq 0$, $(X,Y) \in \mathcal{E}'$,

M4' If $\|x\| \leq X$, $\|y\| \leq Y$, $(X,Y) \in \mathcal{E}'$ then $\|\psi(x,y)\| \leq \Psi(X,Y)$,

are easily deduced as before.

The set of values X for which K5' holds may be characterized by defining the set in $X'_+ \times Y'_+$

$$\Delta = \{(X,Y) \in \mathcal{E}' \mid Y \geq \Psi(X,Y)\}$$

which has projection on to the $X_+$-space

$$\Delta_1 = \{X \in X'_+ \mid \exists Y \in Y'_+ : (X,Y) \in \mathcal{E}' \ \& \ Y \geq \Psi(X,Y)\}.$$

When $X \in \Delta_1$, Kantorovich's condition K5' holds and the following lemma is immediate.

*Lemma 6.1.*

(a) For each $X \in \Delta_1$, the function $F(X) = \Psi(X,Y)$ maps the order interval $\mathcal{D}' = <0,\tilde{Y}>$ into itself.

(b) For such an X, the function $f(y) = \psi(x,y)$ for all values x such that $\|x\| \leq X$, maps the set $\mathcal{D} = \{y \mid \|y\| \in \mathcal{D}'\}$ into itself.

The following iterative processes are now introduced:

$$y^{(0)}(x) = 0$$
$$y^{(p+1)}(x) = \psi(x,y^{(p)}(x)) , \quad p \in \mathbb{N} ,$$

and

$$Y^{(0)}(X) = 0$$
$$Y^{(p+1)}(X) = \Psi(X,Y^{(p)}(X)) , \quad p \in \mathbb{N} .$$

By the last lemma, the iterates are well defined for $\|x\| \leq X \in \Delta_1$.

The following results are then merely a restatement of results in the last section.

*Lemma 6.2.* For $X \in \Delta_1$ the sequence $\{Y^{(p)}\}_{p \in \mathbb{N}}$ is monotonic increasing and bounded.

*Lemma 6.3.* If $\|x\| \le X \in \Delta_1$, the sequences $\{y^{(p)}(x)\}_{p \in \mathbb{N}}$, $\{Y^{(p)}(X)\}_{p \in \mathbb{N}}$ satisfy

(a) $\|y^{(p)}(x)\| \le Y^{(p)}(X)$, $p \in \mathbb{N}$

(b) $\|y^{(p+1)}(x) - y^{(p)}(x)\| \le Y^{(p+1)}(X) - Y^{(p)}(X)$, $p \in \mathbb{N}$.

*Theorem 6.1.* If $\|x\| \le X \in \Delta_1$, the sequences $\{y^{(p)}(x)\}_{p \in \mathbb{N}}$, $\{Y^{(p)}(X)\}_{p \in \mathbb{N}}$ satisfy

(a) $Y^{(p)}(X) \uparrow Y(X)$ where $Y(X) = \Psi(X, Y(X))$

(b) $y^{(p)}(x) \overset{(bk)}{\to} y(x)$ where $y(x) = \psi(x, y(x))$

and $\|y(x)\| \le Y(X)$.

*Definition 6.1.* The solution $Y = Y(X)$ will be called the *principal solution* of the equation $Y = \Psi(X, Y)$.

*Theorem 6.2.* Suppose that for a value $X \in X'_+$, a solution $Y'$ of $Y' = \Psi(X, Y')$ exists. Then the principal solution $Y(X)$ also exists and $Y(X) \le Y'$.

*Proof:* the Kantorovich argument may be applied to the order interval $\langle 0, Y' \rangle$.

Let $\Gamma$ be the graph of $Y = \Psi(X, Y)$ in $X'_+ \times Y'_+$:

$$\Gamma = \{(X,Y) \in \mathcal{E}' \mid Y = \Psi(X,Y)\}$$

and $\Gamma_1$ be its projection on to the $X_+$ space:

$$\Gamma_1 = \{X \in X'_+ \mid \exists Y \in Y_+ : (X,Y) \in \Gamma\}$$

Clearly $\Delta \supset \Gamma$ and $\Gamma$ forms the boundary of $\Delta$. Also $\Delta_1 \supset \Gamma_1$. But conversely, if $X \in \Delta_1$ then by theorem 6.1 also $X \in \Gamma_1$. So $\Delta_1 \supset \Gamma_1$. Consequently $\Delta_1 = \Gamma_1$ i.e. the sets $\Delta, \Gamma$ have the same $X_+$-projection.

The principal solution corresponds to the part

$$\Gamma' = \{(X,Y) \in \Gamma \mid Y = Y(X)\}$$

of the graph $\Gamma$. $\Gamma$ and $\Gamma'$ have the same projection $\Gamma_1$ on the $X'_+$-space. From theorem 6.2 it is seen that $\Gamma'$ is the "lower" part of $\Gamma$ in the sense of partial ordering.

The following property of the principal solution may be noted.

*Lemma 6.4.* $Y(X)$ is isotone increasing.

*Proof:* if $X, X' \in \Delta_1 = \Gamma_1$ and $X \leq X'$ sequences $\{Y^{(p)}(X)\}_{p \in \mathbb{N}}$, $\{Y^{(p)}(X')\}_{p \in \mathbb{N}}$ may be defined. Examining the proof of lemma 6.2. it is seen that $Y^{(p)}(X) \leq Y^{(p)}(X')$, $p \in \mathbb{N}$. Taking limits, $Y(X) \leq Y(X')$.

*Remark:* suppose that $\Psi(X,Y)$ is defined for all $X \geq 0$, $Y \geq 0$, i.e. $\varepsilon = X'_+ \times Y'_+$. Then $Y^{(p)}(X) \uparrow +\infty$ as $p \to \infty$ when $X \in X'_+ \setminus \Delta_1$ ($= X'_+ \setminus \Gamma_1$) since the iterates are then well-defined but the sequence cannot converge otherwise the point $X$ would be the projection of a point on $\Gamma$, i.e. it would be true that $X \in \Gamma_1$.

## 7. Generalization of Hille's Theorem

In this sections will be considered equations

$$Y = \Psi(X,Y) \qquad y = \psi(x,y)$$

where $\Psi, \psi$ are subject to conditions (a), (b), (c) below.

(a) $\qquad \Psi(X,Y) = (o) - \sum\sum_{m+n\geq 1} \frac{1}{m!n!} A_{mn}(X,\ldots,X; Y,\ldots,Y)$

$A_{mn}$ multilinear, completely symmetrical in X and Y variables, continuous for monotone convergence,

$$A_{01} = 0,$$

Series for $\Psi(X,Y)$ (o)-convergent for $(X,Y) \in \mathcal{E}'$ where $\mathcal{E}'$ is an (o)-star in $X'_+ \times Y'_+$.

(b) $\qquad \psi(x,y) = (bk) - \sum\sum_{m+n\geq 1} \frac{1}{m!n!} a_{mn}(x,\ldots,x; y,\ldots,y)$

$a_{mn}$ multilinear, completely symmetrical in x and y variables,

$$a_{01} = 0,$$

(c) $\qquad \psi \ll \Psi$.

From these conditions it follows that the series for $\psi$ is (bk)-convergent when $(\|x\|,\|y\|) \leq (X,Y) \in \mathcal{E}'$ i.e. on the set $\mathcal{E} = \{(x,y) \mid (\|x\|,\|y\|) \in \mathcal{E}'\}$

The formal series solutions of the equations will, as in section 4, be denoted respectively by

$$Y = \sum_{n=1}^{\infty} \frac{1}{n!} B_n(X,\ldots,X) = \Phi(X)$$

$$y = \sum_{n=1}^{\infty} \frac{1}{n!} b_n(x,\ldots,x) = \varphi(x)$$

From lemma 4.2, $\Phi \gg \varphi$.

The conditions K1'-K4', M1'-M3' of section 6 are easily verified and so the results of that section may be used.

*Lemma 7.1.* Under the above conditions

(a) The $\{Y^{(p)}(X)\}_{p\in\mathbb{N}}$ form a sequence of (o)-analytic functions with power series of positive type with region of convergence $X \in \Delta_1 \;(= \Gamma_1)$.

(b) The iterates $\{y^{(p)}(x)\}_{p\in\mathbb{N}}$ form a sequence of (bk)-analytic functions with region of convergence $\{x \mid \|x\| \le X \in \Delta_1 \;(= \Gamma_1)\}$

*Proof:*

(a) Induction. $Y^{(0)}(X) = 0$ is certainly (o)-analytic. Assuming $Y^{(p)}(X)$ (o)-analytic in $X \in \Delta_1 = \Gamma_1$, it satisfies by theorem 6.1 the inequality $0 \le Y^{(p)}(X) \le Y(X)$ on $X \in \Delta_1 = \Gamma_1$ and so (theorem 4. ) $Y^{(p+1)}(X) = \Psi(X, Y^{(p)}(X))$ is also (o)-analytic on $X \in \Delta_1 = \Gamma_1$ since then $(X, Y(X)) \in \mathcal{E}'$.

(b) This part follows immediately using the same argument and lemma 4.3.

*Lemma 7.2.* Under the same conditions, the power series expansions of $Y^{(p)}(X)$, $y^{(p)}(x)$ coincide with the first p terms of the formal series solutions, i.e.

$$Y^{(p)}(X) = \sum_{n=1}^{p} \frac{1}{n!} B_n(X,\ldots,X) + \text{terms of positive type of degree} > p$$

$$y^{(p)}(x) = \sum_{n=1}^{p} \frac{1}{n!} b_n(x,\ldots,x) + \text{terms of degree} > p$$

*Proof:* straightforward induction using the formulae of section 4.

*Theorem 7.1.* (Generalized Hille Theorem)
Considering the analytic implicit function equations $y = \psi(x,y)$ and $Y = \Psi(X,Y)$ defined at the beginning of this section,

(a) The solution $Y = \Phi(X)$ of $Y = \Psi(X,Y)$ found by series substitution is (o)-convergent for $X \in \Delta_1 \;(= \Gamma_1)$ and, for these values of $X$, coincides with the principal solution $Y(X)$.

(b) The solution $y = \varphi(x)$ of $y = \psi(x,y)$ found by series substitution has Cauchy's majorant $\Phi$ and is (bk)-convergent when $\|x\| \le X \in \Delta_1 \;(= \Gamma_1)$.

*Proof:*

(a) For $X \in \Delta_1 \;(= \Gamma_1)$ the principal solution $Y(X)$ is defined and, by lemma 7.2,

$$\sum_{n=1}^{p} \frac{1}{n!} B_n(X,\ldots,X) \le Y^{(p)}(X) \le Y(X)$$

This is true for all integers $p \geq 1$. Hence

$$\Phi(X) = (o) - \sum_{n=1}^{\infty} \frac{1}{n!} B_n(X,\ldots,X)$$

is (o)-convergent in a convergence region $X \in \Delta_1 (= \Gamma_1)$. It must satisfy

$$\Phi(X) \leq Y(X)$$

Since $(X,Y(X)) \in \mathcal{E}'$ for all $X \in \Delta_1 (= \Gamma_1)$ it follows (lemma 4.3) that the series for $Y = \Phi(X)$ may be substituted in the series for $\Psi(X,Y)$ giving an (o)-analytic function $\Psi(X,\Phi(X))$. Since, formally, $\Phi(X) = \Psi(X,\Phi(X))$, it follows that this equation is in fact satisfied and $Y = \Phi(X)$ gives a solution of $Y = \Psi(X,Y)$. Now by the minimal property of the principal solution of the principal solution $Y(X)$ (theorem 6.1) it follows that $\Phi(X) \leq Y(X)$ for $X \in \Delta_1 (= \Gamma_1)$ and since the reverse inequality has already been shown, it must be that $\Phi(X) = Y(X)$ for $X \in \Delta_1 (= \Gamma_1)$.

(b)  The proof of this part follows immediately.

*Remark*: if the series for $\Psi(X,Y)$ is (o)-convergence for all $X \geq 0$, $Y \geq 0$ (i.e. $\mathcal{E}' = X'_+ \times Y'_+$) then the series solution $Y = \Phi(X)$ is divergent when $X \in X'_+ \setminus \Delta_1 (= X'_+ \setminus \Gamma_1)$ cf. the remark at the end of the previous section. The sequence $Y^{(p)}(X)$ is well-defined but cannot converge.

Acknowledgement

The author expresses his appreciation to Prof.dr.ir. M.L.J. Hautus for the possibility of working in the Mathematics Department of TH Eindhoven where this report was completed in connexion with a contract to investigate the use of Riesz space in control theory.

The author also thanks the typing staff for their willing help in the final production of the report.

## References

1. Cristescu R.   Ordered Vector Spaces and Linear Operators
   English edition: Tunbridge Wells, (G.B.), 1976 (Abacus Press)

2. Goursat E.   A Course in Mathematical Analysis
   English edition: Boston 1904 (Ginn)   Reprint: New York 1950 (Dover)

3. Hille E.   Functional Analysis and Semigroups
   Amer. Math. Soc. Colloq. Pub. 1948

4. Hille E.   Analytic Function Theory
   Boston etc. 1959 (Ginn)   See theorem 9.4.6. on p. 274 of vol I

5. Kantorovich L.V.   The method of successive approximation for functional equations.
   Acta Math. 1939, $\underline{71}$, 63-97

6. Kantorovich L.V., Vulikh B.Z. & Pinsker A.G.   Partially ordered groups and partially ordered linear spaces.
   Usp. Mat. Nauk. 1951, $\underline{6}$, no. 3 (43), 31-98
   Amer. Math. Soc. Trans., Sept. 1963, vol. 27, ser 2, 51-124

7. Mazur S. & Orlicz V.   Fundamental properties of polynomial operations.
   Studia Math. 1935, vol 5, 60-68

8. Michal A.D.   Differential Calculus in Banach Spaces.
   Paris 1958 (Gauthier-Villars)   (French)

9. Riordan J.   Combinatorial Identities.
   New York 1968 (Wiley)

10. Vulikh B.Z.   Introduction to the Theory of Partially Ordered Spaces.
    English edition: Groningen 1967 (Wolters-Noordhoff)